\documentclass[11pt,amsfonts]{article}
\usepackage[utf8]{inputenc}
\usepackage[T1]{fontenc} 
\usepackage{tikz}
\usepackage[english]{babel} 
\usepackage{mathrsfs}
\usepackage{amssymb}
\usepackage[cmex10]{amsmath}
\usepackage{amsthm}
\usepackage{wrapfig}
\usepackage{latexsym}
\usepackage{amsfonts}
\usepackage{color}
\usepackage[font=small, labelfont=bf]{caption}
\usepackage{ctable}
\usepackage{floatrow}
\usepackage[colorlinks=true, citecolor=red, linkcolor=blue]{hyperref}
\usepackage{graphicx}
\usepackage{subfig}
\usepackage{verbatim}
\usepackage{epstopdf}
\usepackage{float}

\usepackage[top=2cm , bottom=2cm , left=4cm ,right=4cm ]{geometry}

\usepackage{cases}

\theoremstyle{plain}
\newtheorem{definition}{Definition}
\newtheorem{theorem}{Theorem}

\newtheorem{proposition}{Proposition}

\newtheorem{lemma}{Lemma}
\newtheorem*{pve}{Proof}

\newtheorem{remark}{Remark}

\begin{document}
		\renewcommand{\thefootnote}{\arabic{footnote}}
	\begin{center}
		{\Large \textbf{On Skorokhod Problem  with Two RCLL Reflecting Completely Separated Barriers\footnote{This work is supported by Hassan II Academy of Sciences and Technology }}} \\[0pt]
		~\\[0pt]
		Rachid Belfadli\footnote{Department of Mathematics, Faculty of Sciences and Technology Cadi Ayyad University, Gu{\'e}liz
			40 000 Marrakech Morocco. Email: \texttt{r.belfadli@uca.ma},} Imane Jarni\footnote{Mathematics Department, Faculty of Sciences Semalalia, Cadi Ayyad University, Boulevard Prince Moulay Abdellah,	P. O. Box 2390, Marrakesh 40000, Morocco. 
			E-mail: \texttt{jarni.imane@gmail.com}} and Youssef Ouknine \footnote{Complex Systems Engineering and Human Systems, Mohammed VI Polytechnic University, Lot 660, Hay Moulay Rachid, Ben Guerir 43150, Morocco
				\\
				Mathematics Department, Faculty of Sciences Semalalia, Cadi Ayyad University, Boulevard Prince Moulay Abdellah, P. O. Box 2390, Marrakesh 40000, Morocco. E-mail: \texttt{youssef.ouknine@um6p.ma, ouknine@uca.ac.ma}}%
			\\[0pt]
		
		\textit{Cadi Ayyad University and Mohammed VI Polytechnic University  }\\[0pt]
		~\\[0pt]
	\end{center}

	\begin{abstract}
		In this paper we deal with Skorokhod problem for right continuous left limited (rcll) barriers. We prove existence and uniqueness of the solution when the barriers are only supposed to be rcll and completely separated. Then, we apply our results to prove  existence and uniqueness of the solution of a reflected stochastic differential equation (SDE).
		
	\end{abstract}
	
	{\bf Keyword}: Skorokhod problem; Reflection; Reflecting Stochastic Differential Equations; Approximation method .\\
	\textbf{AMS Subject Classification:} Primary 60K99, Secondary 60H20
	
	\section{Introduction}	
	
	\quad \quad  Equations with reflecting   condition on the boundary of a convex domain, in particular an interval, have been studied by many authors. One barrier reflected equations were introduced firstly by Skorokhod in \cite{12}. He  constructed the solution of SDE staying in the half-line $\mathbb{R}^{+}$ and verifying a reflecting boundary condition at the barrier $0$. Moreover, Chaleyat-Maurel and El Karoui, considered a deterministic version of the reflection problem for continuous functions in \cite{2}. While in \cite{3}, Chaleyat-Maurel, El Karoui and Marchal have extended and proved the existence and uniqueness of the solution of the  reflection problem in $\mathbb{R}^{+}$ for rcll functions. Let us also mention the recent  work   by Slomi{\'n}ski and Wojciechowski \cite{13}, in which the authors have set up and solved a new version of the Skorokhod problem generalizing the one stated in \cite{3}. Indeed, the authors considered a parameter of reflection $(a_{t})_{t\geq 0}$ taking values in $[0,\frac{1}{2}]$; the case  $a=\frac{1}{2}$ corresponds to the classical one given in \cite{3}. As a byproduct,  the authors applied their results to solve reflected SDEs. 
	
	\quad The equations with two reflecting barriers, called also  two sided reflection equations,  have received a lot of attention by many researchers. The case of constants barriers was studied in (\cite{14}, \cite{9}, \cite{6}) while the one of continuous barriers was considered in \cite{8}.
	In addition, Pihlsg\aa{}rd and Glynn in \cite{11} dealt with the case where the barriers are supposed  to be rcll semimartingales $L$ and $U$ satisfying the following hypothesis:  $$(\mathbf{H_{0}}):\;\;\; \text{there exists}\; \epsilon>0 \; \text{such that} \;  U_{t}-L_{t}>\varepsilon. \;\;\;\;\;\;\;\; \text{ for all}\;\; t\in [0,+\infty[.\;\;\;\;\;\;\;\;\;\;\;\;\; $$                                  
	
	In \cite{11}, The authors prove existence and uniqueness of the solution of Skorokhod problem under the assumption $\mathbf{(H_{0})}$. Our main motivation in writing  this paper is to deal with  Skorokhod problem, in finite horizon $[0,T]$, under a rather weaker condition than $(\mathbf{H_{0}})$. More precisely, we solve the reflection Skorokhod problem with respect two  rcll and adapted processes $L$ and $U$, under  the following relaxed  hypothesis:
	
	$$(\mathbf{H_{1}}):\;\;\;\;\; \mathbb{P} \text{-a.s for all}\;\; t\in[0,T],\;\; L_{t}<U_{t} \;\; \text{and}\;\; L_{t^{-}}<U_{t^{-}}. \;\;\;\;\;\  $$
	
	This condition means that the barriers are completely separated. Clearly the hypothesis $(\mathbf{H_{0}})$  implies the hypothesis $(\mathbf{H_{1}})$, but the converse is not true  since $\varepsilon$ will depend on $\omega$. The hypothesis $(\mathbf{H_{1}})$ is inspired from the work of Hamad{\`e}ne, Hassani and Ouknine in \cite{5}, where the hypothesis $(\mathbf{H_{1}})$ is used instead of the so called Mokobodski's condition to solve reflected BSDEs .
	
	As an application of our result and based on  a  Picard-type approximation  we establish the existence and uniqueness of the solution of a class of reflected stochastic differential equations. Namely,  we consider the following  type of  reflecting equations with respect to two rcll  barriers $L$ and $U$:
	$$ X_{t}=H_{t}+\int_{0}^{t}\sigma(s,X_{s})dB_{s}+\int_{0}^{t}a(s,X_{s})ds+K_{t}^{+}-K^{-}_{t},\;\;\;  $$
	
	where $(B_{t})_{0\leq t\leq T}$ is a standard Brownian motion, the process $(H_{t})_{0\leq t\leq T}$ is a  rcll semimartingale with $L_{0}\leq H_{0}\leq U_{0}$ and $K=K^{+}-K^{-}$ is the rcll minimal process with bounded variation making $X$ in the interval  $[L_{t},U_{t}]$.
	
	This paper is organized as follow : In  Section 2, we set the Skorokhod problem we deal with, and recall the related results of existence and uniqueness when the barriers satisfy the condition $\mathbf{(H_{0})}$. Then we state and prove our first main result, namely Theorem 1, showing the existence and uniqueness of the solution of the Skorokhod problem  under the relaxed condition $\mathbf{(H_{1})}$. In Section 3 we apply our result to study existence and uniqueness for a class of SDEs with two reflecting barriers.  Finally, we present in an Appendix a slight generalization of the result in [\cite{8},p 270-271 ] to the case when the barriers are only rcll adapted processes. That is,  we will show  the existence and uniqueness of the solution of the Skorokhod problem under the assumption $\mathbf{(H_{0})}$ and when the barriers are only rcll.\\
	
	\large{\textbf{ Notations.}}\\
	
	Throughout this section, we are given a probability space $(\Omega,\mathcal{F},\mathbb{P}) $  equipped with a complete filtration $(\mathcal{F}_{t})_{t\in[0,T]}$ that we suppose satisfying  the usual conditions (completion and right continuity).
	
	\begin{itemize}
		\item[$\bullet$] For a rcll $X$, we denote  $X_{t^{-}}:=\lim\limits_{s\nearrow t} X_{s}$ its left limit, $\Delta X_{t}:=X_{t}-X_{t^{-}}$ and $\mathbb{E}_{s}(X_{t}):=\mathbb{E}(X_{t}/ \mathcal{F}_{s}) $ for $s\leq t$. For any $x\in\mathbb{R}$, $x^{+}:=\max(x,0)$.
		\item[$\bullet$] The total variation of a process $K$ with bounded variation on $[0,T]$ will be denoted by $\text{Var}_{[0,T]}(K)$.
		\item[$\bullet$] $\mathbb{L}^{2}$: The set of real valued $\mathcal{F}$-measurable square integrable random variables. 
		\item[$\bullet$] $\mathcal{S}^{2}$: The complete space of real valued   and $\mathcal{F}_t$-adapted rcll processes $(X_{t})_{0\leq t\leq T}$ such that   $||X||_{\mathcal{S}^{2}}=|| \sup_{0\leq t\leq T} |X_{s}|\; ||_{\mathbb{L}^{2}}<+\infty$. 
	\end{itemize} 
	
	\section{Skorokhod problem with two rcll reflecting separately barriers}
	Let us first give the definition of the reflected Skorokhod problem between two barriers. We are given three rcll adapted processes $Y$, $L$ and $U$ such that $L\leq U$ and $L_{0}\leq Y_{0}\leq U_{0}$. 
	\begin{definition}\label{def1}
		A triplet of processes $(X, K^{+}, K^{-})$ is said to be a solution of the Skorokhod problem associated with $Y$ and reflected at the barriers $L$ and $U$, that we shortly labeled  $SPR(Y, L, U)$, if the processes $K^{+}$ and $K^{-}$ are rcll adapted, nonnegative, and  increasing such that $K_{0}^{+}=K_{0}^{-}=0$ and for all $t\in [0,T]$ the following holds:
		\begin{equation} X_{t}=Y_{t}+K^{+}_{t}-K^{-}_{t} \label{1}, \end{equation}
		\begin{equation} L_{t} \leq X_{t}\leq U_{t}, \label{2} \end{equation}
		and
		\begin{equation}\displaystyle\int_{0}^{T} (X_{s}-L_{s})\; \rm{d}K_{s}^{+} =\displaystyle\int_{0}^{T} (U_{s}-X_{s}) \; \rm{d}K_{s}^{-}=0.\label{3} \end{equation}
		
		Here the integrals are taken to be  understood in Stieltjes's sens with respect to the increasing processes $K^{+}$ and $K^{-}$. 
	\end{definition} 
	
	Notice that in the Definition \ref{def1} we do not require the property of being a semimartingale neither for $Y$, $L$ nor for $U$. In fact, and in contrast of \cite{11},  the following result borrowed from \cite{11},  shows that if the processes $Y$, $L$ and $U$ are assumed to be  semimartingales and if the assumption $\mathbf{(H_{0})}$ holds, then the Skorokhod problem $SPR(Y,L,U)$ has a unique solution.

	\begin{proposition}[{\cite{11}}]\label{pro1}
		If the processes are rcll semimartingales and if the assumption $(\mathbf{H_{0}})$ is fulfilled, then
		the Skorokhod problem  $SPR(Y,L,U)$ has a unique solution.
	\end{proposition}
	
	However, following the lines of the proof of Proposition 2 in {\cite{11}}, it is seen  that the semimartingale assumption is no more needed neither on the process $Y$ nor on the barriers $L$ and $U$,  and  then the conclusion of the proposition still valid in this case. For the convenience of the reader we give a complete proof in the Appendix.\\
	
	As pointed out in the introduction our main objective, and  this is  the novelty of this paper, is to show that the Skorokhod problem $SPR(Y, L, U)$ has a unique solution under the weaker condition $\mathbf{(H_{1})}$ compared to $\mathbf{(H_{0})}$.\\
	
	We begin with the following lemma which will be used in the proof of our main result.
	
	\begin{lemma}\label{lem1}
		Let $Y$, $L$, $U$  and $\tilde{U}$ be rcll adapted processes such that $L\leq U$, $L\leq \tilde{U}$, $L_{0}\leq Y_{0}\leq U_{0}$ and $L_{0}\leq Y_{0}\leq \tilde{U}_{0}$. If  $(X, K^{+}, K^{-})$ \big(resp. $(\tilde{X},\tilde{K}^{+},\tilde{K}^{-})$ \big) is a solution of the Skorokhod problem $SPR(Y, L, U)$ \big(resp. $SPR(Y, L, \tilde{U})$ \big), then we have, $\mathbb{P}$-a.s for all $t\in [0,T]$,
		\begin{equation}
		(X_{t}-\tilde{X}_{t})^{2}\leq 2\displaystyle\int_{0}^{t}(\tilde{U}_{s}-U_{s}) \rm{d}K^{-}_{s}+ 2\displaystyle\int_{0}^{t}(U_{s}-\tilde{U}_{s}) \rm{d}\tilde{K}^{-}_{s}\;\; \mathbb{P}\text{-a.s} .  \label{4} 
		\end{equation}
	\end{lemma}
	\begin{pve} From (\ref{1}) we have  $X-\tilde{X}=K^{+}-K^{-}-\tilde{K}^{+}+\tilde{K}^{-}$, which is a process of bounded variation. Then we apply the integration by part formula  we get:
		$$
		\begin{array}{lll}
		(X_{t}-\tilde{X}_{t})^{2}&=&2 \displaystyle\int_{0}^{t} (X_{s}-\tilde{X}_{s})\rm{d}(X_{s}-\tilde{X}_{s}) - \displaystyle\sum_{0< s \leq t}
		\Delta (X-\tilde{X})_{s}^{2}\\
		&\leq& 2\displaystyle\int_{0}^{t} (X_{s}-\tilde{X}_{s})\rm{d}K_{s}^{+}-2\displaystyle\int_{0}^{t} (X_{s}-\tilde{X}_{s})\rm{d}K^{-}_{s}-2\displaystyle\int_{0}^{t} (X_{s}-\tilde{X}_{s})\rm{d}\tilde{K}^{+}_{s}\\
		& & +2\displaystyle\int_{0}^{t} (X_{s}-\tilde{X}_{s})\rm{d}\tilde{K}^{-}_{s}\\
		
		&=& 2\displaystyle\int_{0}^{t} (X_{s}-L_{s})\rm{d}K_{s}^{+} +
		2\displaystyle\int_{0}^{t} (U_{s}-X_{s})\rm{d}K_{s}^{-}+
		2\displaystyle\int_{0}^{t} (\tilde{X}_{s}-L_{s})\rm{d}\tilde{K}^{+}_{s}\\
		& &+	2\displaystyle\int_{0}^{t}(\tilde{U}_{s}-\tilde{X}_{s}) \rm{d}\tilde{K}^{-}_{s}+ 
		2\displaystyle\int_{0}^{t} (L_{s}-X_{s})\rm{d}\tilde{K}^{+}_{s}+
		2\displaystyle\int_{0}^{t} (L_{s}-\tilde{X}_{s})\rm{d}K_{s}^{+}\\ 
		& &+		2\displaystyle\int_{0}^{t} (\tilde{X}_{s}-U_{s})\rm{d}K_{s}^{-}
		+
		2\displaystyle\int_{0}^{t}(X_{s}-\tilde{U}_{s}) \rm{d}\tilde{K}^{-}_{s}.
		\\
		\end{array} $$
		But since the first four integrals are equal to zero thanks to (\ref{3}), and the fifth and the sixth integrals are less than $0$ thanks to (\ref{2}), we obtain:
		$$\begin{array}{lll}
		(X_{t}-\tilde{X}_{t})^{2}&\leq& 2\displaystyle\int_{0}^{t} (\tilde{X}_{s}-U_{s})\rm{d}K_{s}^{-}
		+ 2\displaystyle\int_{0}^{t}(X_{s}-\tilde{U}_{s}) \rm{d}\tilde{K}^{-}_{s}\\
		& \leq & 2\displaystyle\int_{0}^{t}(\tilde{U}_{s}-U_{s}) \rm{d}K^{-}_{s}+ 2\displaystyle\int_{0}^{t}(U_{s}-\tilde{U}_{s}) \rm{d}\tilde{K}^{-}_{s},
		\end{array}$$
		
		where we have used the (\ref{2}) in the last inequality.
		This completes the proof of Lemma 1 .
	\end{pve}
	\begin{remark}
		By a symmetric argument one can check  that if  $(X, K^{+}, K^{-})$ \big(resp. $(\tilde{X},\tilde{K}^{+},\tilde{K}^{-})$ \big) is a solution of the Skorokhod problem $SPR(Y, L, U)$ \big(resp. $SPR(Y, \tilde{L}, U)$ \big), then we have, $\mathbb{P}$-a.s for all $t\in [0,T]$,
		\begin{equation*}
		(X_{t}-\tilde{X}_{t})^{2}\leq 2\displaystyle\int_{0}^{t}(L_{s}-\tilde{L}_{s}) \rm{d}K^{+}_{s}+ 2\displaystyle\int_{0}^{t}(\tilde{L}_{s}-L_{s}) \rm{d}\tilde{K}^{+}_{s}\;\; \mathbb{P}\text{-a.s} .  
		\end{equation*}
	\end{remark}
	We are now in position to state the main result of this Section.
	
	\begin{theorem}\label{Th1}
		Let $Y$, $L$, and $U$ be rcll adapted processes such that $L_{0}\leq Y_{0}\leq U_{0}$ and $L$ and $U$ satisfying the assumption $(\mathbf{H_{1}})$. Then, the Skorokhod problem  $SPR(Y, L, U)$ has a unique solution.
	\end{theorem}
	\begin{pve}
		\textit{(i)-First we prove the existence part}.\\  We set for any integer $n\geq 1$,
		$U^{n}:=\max(U,L+\frac{1}{n})$. Since 
		$U^{n}-L\geq \frac{1}{n}$ on $[0,T]$, the result of Proposition \ref{proApp} in the Appendix  may be applied and then the Skorokhod problem $SPR(Y, L, U^{n})$  has a unique solution $(X^{n},K^{n,+},K^{n,-})$. In order to show that the sequence $(X^{n}, K^{n,+}, K^{n,-})$ converges to the solution of our Skorokhod problem  $SPR(Y, L, U)$, we first prove  that:
		\begin{equation}
		\mathbb{P}\Big[\bigcup_{N\geq 1}\bigcap_{n\geq N} \big\{\sup_{0\leq t \leq T}(U_{t}^{n}-U_{t})=0 \big\}\Big]=1. \label{5}
		\end{equation}
		This claim means that the sequence $\sup_{0\leq t \leq T} (U_{s}^{n}-U_{s})$ is of stationary type. \\
		Now, since $(U^{n})_{n\geq 1}$ is decreasing, we only need to prove that $\mathbb{P}(E)=0 $, where the set $E$ is given by:
		$$E=\{\omega \; : \:\; \forall n\geq 1\;\; \exists t\in [0,T]  \;\; U_{t}^{n}(\omega)>U_{t}(\omega)\;\}.$$
		Let $\omega \in E$ and $ n\geq 1 $. Then there exists a sequence  $(t_{n})_{n\geq 0}$ in $[0,T]$ such that    $U_{t_{n}}^{n}(\omega)>U_{t_{n}}(\omega)$. Thus
		$U_{t_{n}}^{n}(\omega)=L_{t_{n}}(\omega)+\frac{1}{n}$
		and so  $0 \leq U_{t_{n}}(\omega)-L_{t_{n}}(\omega)\leq \frac{1}{n}$. Therefore,
		\begin{equation}
		\lim\limits_{n\rightarrow +\infty} [U_{t_{n}}(\omega)-L_{t_{n}}(\omega)]=0. \label{ 6}
		\end{equation}
		On the other hand, from the sequence $(t_{n})$ in $[0,T]$, we can subtract a subsequence $(t_{n_{k}})_{k\geq 1}$ converging to some element $t$ in $[0,T]$. Now, since either one of the  two sets $\{k\in\mathbb{N}: t_{n_{k}}>t\}$ and $\{k\in\mathbb{N}: t_{n_{k}}\leq t  \}$ is at least infinite, we deduce that
		
		\begin{equation}
		\begin{split}
	&	\lim\limits_{k\rightarrow +\infty} [U_{t_{n_{k}}}(\omega)-L_{t_{n_{k}}}(\omega)]=U_{t}(\omega)-L_{t}(\omega)\\
	&	\; \text{or} \;\\
	&	\lim\limits_{k\rightarrow +\infty} [U_{t_{n_{k}}}(\omega)-L_{t_{n_{k}}}(\omega)]=U_{t^{-}}(\omega)-L_{t^{-}}(\omega). \label{7}
		\end{split}
		\end{equation}
		
		It follows from (\ref{ 6}) and (\ref{7}) that:
		$$U_{t}(\omega)=L_{t}(\omega) \;\;\text{or}\;\; U_{t^{-}}(\omega)=L_{t^{-}}(\omega).$$
		Whence $$ E\subseteq \{ \omega \;:\; \exists t\in [0,T] \; U_{t}=L_{t} \;\text{ou}\; U_{t^{-}}=L_{t^{-}} \} , $$ which is $\mathbb{P}$-null set. because of the assumption $(\mathbf{H_{1}})$. \\
		
		Let us now focus on the convergence of the sequence $(X^{n},K^{n,+},K^{n,-})$. By \eqref{5}, we have 
		$\mathbb{P}\text{-a.s}$,
		\begin{equation}
		\exists N\geq 1 \;\; \forall  n\geq N  \;\; \sup_{0\leq t \leq T}(U_{t}^{n}-U_{t})=0. \label{8}
		\end{equation}

		Therefore, using Lemma 1, we get for all $n,m\geq N$
		\begin{equation} 
		\sup_{0\leq t \leq T}(|X_{t}^{n}-X_{t}^{m}|)=0.\label{9}
		\end{equation}	
		That is the sequence $X^{n}$ is of stationary type and then converging in particular to a certain limit $X$ which is rcll and satisfying $ L\leq X \leq U$ due to (\ref{2}).\\
		Now going back to (\ref{8}), we have for all $n,m\geq N$ and $t\in[0,T]$, $U_{t}^{n}=U_{t}^{m}=U_{t}$, which implies together with (\ref{9}) that $K^{n,+}_{t}-K^{m,+}_{t}=K^{n,-}_{t}-K^{m,-}_{t}$ and therefore:
		$$\begin{array}{lll}
		& &\displaystyle\int_{0}^{T}(U_{s}-L_{s})\rm{d}(K^{n,+}_{s}-K^{m,+}_{s})\\
		&=&\displaystyle\int_{0}^{T}(U_{s}-X_{s})\rm{d}(K^{n,+}_{s}-K^{m,+}_{s})+\displaystyle\int_{0}^{T}(X_{s}-L_{s})\rm{d}(K^{n,+}_{s}-K^{m,+}_{s})\\
		&=&\displaystyle\int_{0}^{T}(U_{s}-X_{s})\rm{d}(K^{n,-}_{s}-K^{m,-}_{s})+\displaystyle\int_{0}^{T}(X_{s}-L_{s})\rm{d}(K^{n,+}_{s}-K^{m,+}_{s})\\
		&=&\displaystyle\int_{0}^{T}(U_{s}^{n}-X_{s}^{n})\rm{d}K^{n,-}_{s}-\displaystyle\int_{0}^{T}(U_{s}^{m}-X_{s}^{m})\rm{d}K^{m,-}_{s}\\
		&  &  +\displaystyle\int_{0}^{T}(X_{s}^{n}-L_{s})\rm{d}K^{n,+}_{s} + \displaystyle\int_{0}^{T}(X_{s}^{m}-L_{s})\rm{d}K^{m,+}_{s}=0,
		\end{array}
		$$
		which implies that $K^{n,+}=K^{m,+}$ and $K^{n,-}=K^{m,-}$ since the barriers are completely separated.\\
		We denote $K^{+}:=\lim\limits_{n\rightarrow+\infty} K^{n,+}$ and $K^{-}:=\lim\limits_{n\rightarrow+\infty} K^{n,-}$.
		so $K^{+}$ and $K^{-}$ are increasing, rcll and adapted  processes. On the other hand, for $n\geq N$ we have:
		
		$$\displaystyle\int_{0}^{T}(X_{s}-L_{s})\rm{d}K^{+}_{s}=\displaystyle\int_{0}^{T}(X_{s}^{n}-L_{s})\rm{d}K^{n,+}_{s}=0,$$
		and
		$$\displaystyle\int_{0}^{T}(X_{s}-U_{s})\rm{d}K^{-}_{s}=\displaystyle\int_{0}^{T}(X_{s}^{n}-U_{s})\rm{d}K^{n,-}_{s}=0.$$
		
		Therefore  $K^{+}$ and $K^{-}$ satisfy $\eqref{3}$ and finally $(X,K^{+},K^{-})$ is a solution of $SPR(Y, L, U)$.
		\\
		\textit{(ii)}- For the uniqueness part, if $(X,K^{+},K^{-})$ and $(X^{*},K^{+,*},K^{-,*})$ are two solutions associated to $SPR(Y, L, U)$, then by Lemma 1, for all $t\in [0,T]$  we have $ X_{t}=X^{*}_{t} $. But since $X$ and $X^{*}$ are rcll, they are indistinguishable.
		and On the other hand we have: 
		\begin{equation*}
		\begin{split}
		\int_{0}^{T}(U_{s}-L_{s})d(\rm K^{+}_{s}-K_{s}^{+,*})&=\int_{0}^{T}(U_{s}-X_{s})d(\rm K^{+}_{s}-K_{s}^{+,*})\\
		&=\int_{0}^{T}(U_{s}-X_{s})d(\rm K^{-}_{s}-K_{s}^{-,*})=0
				\end{split}
		\end{equation*}
		Since $L$ and $U$ are completely separated then $K^{+}=K^{+,*}$ and	$K^{-}=K^{-,*}$.

		Which completes the proof of the Theorem.

	\end{pve}
	
	\begin{remark}
		Observe that by using  Remark 1 it is possible to take the barriers $(L^{n}, U)$ instead of $(L,U^{n})$, where $L^{n}=\min(L,U-\frac{1}{n})$ in the above argument. The proof follows similarly.
	\end{remark}
	We close this section by giving an equivalent formulation of the assumption $\mathbf{(H_{1})}$.	
	\begin{proposition}
		If $L$ and $U$ are two rcll adapted processes such that\\ $L_{0}=L_{0^{-}}$ and $U_{0}=U_{0^{-}}$, then the hypothesis $\mathbf{(H_{1})}$ is equivalent to \\$ \inf_{t\in [0,T]}(U_{t}-L_{t})>0$, $\mathbb{P}$-a.s.
	\end{proposition}
	\begin{pve}
		Since the sufficiency follows immediately, we only show that :\\ $ \inf_{t\in [0,T]}(U_{t}-L_{t})>0,$ holds $\mathbb{P}\text{-a.s.}$ whenever $\mathbf{(H_{1})}$ is satisfied. By the sequential characterization of the infimum, we can find a sequence $(s_{n})_{n\geq 0}$ in $[0,T]$ such that: $$\lim\limits_{
			n\rightarrow+\infty}(U_{s_{n}}-L_{s_{n}})=\inf_{t\in [0,T]}(U_{t}-L_{t}).$$ But there exists a subsequence $(s_{n_{k}})_{k\geq 0}$  converging to some $s\in [0,T]$, and again, as in the proof of Theorem 1, we have $$\lim\limits_{k\rightarrow+\infty}(U_{s_{n_{k}}}-L_{s_{n_{k}}})=  U_{s}-L_{s}\;\text{or}\; \lim\limits_{k\rightarrow+\infty}(U_{s_{n_{k}}}-L_{s_{n_{k}}})=  U_{s^{-}}-L_{s^{-}},$$
		which implies that 
		$ \inf_{t\in [0,T]}(U_{t}-L_{t})>0$, $\mathbb{P}$-a.s.

	\end{pve}
	
	\section{SDEs with two reflecting barriers }
	In in this section, we consider a complete probability space endowed with a standard Brownian motion  $B=(B_{t})_{0\leq t\leq T }$ and the usual augmented filtration of $B$, $\{\mathcal{F}_{t}\}_{0\leq t\leq T}$ .

	We consider the following norm introduced in \cite{10}:
	$$||X||^{2}:=||X||_{\mathcal{S}^{2}}^{2}+\displaystyle\sup_{\pi}\mathbb{E}\Big[\Big( \displaystyle\sum_{i=0}^{m-1}|\mathbb{E}_{\tau_{i}}(X_{\tau_{i+1}})-X_{\tau_{i}}|\Big)^{2}\Big],$$
	where $X$ is a progressively measurable rcll process, and the supremum is taken over all partitions $\pi:0=\tau_{0}\leq...\leq\tau_{m}=T$ of $[0,T]$, for some stopping times $\tau_{0},...,\tau_{m}$. Moreover, if $X$ is a rcll semimartingale, then $X$ has the following decomposition: $$X_{t}=X_{0}+M_{t}+A_{t} ,$$ where $M$ is a local martingale and $A$ is a process of finite variation.\\
	
	The following result is borrowed from \cite{10} gives some characterizing results for semimartingales $X$ satisfying:  
	$$ \mathbb{E}\big[|X_{0}|^{2}+<M>_{T}+[\text{Var}_{[0,T]}(A)]^{2}\big] <+\infty.$$
	
	\begin{theorem}[{\cite{10}}]
		Let $X$ be a rcll adapted process. Then, the following assertions hold:
		\begin{itemize}
			\item[(i)] If $X$ is a semimartingale  such that $$ \mathbb{E}\big[|X_{0}|^{2}+ < M >_{T}+[\text{Var}_{[0,T]}(A)]^{2}\big] <+\infty,$$ then there exists universal constants  $0<c< C$ such that:
			$$ c||X||\leq \mathbb{E}\big[|X_{0}|^{2}+ <M>_{T}+[\text{Var}_{[0,T]}(A)]^{2}\big]\leq C||X||. $$
			
			\item[(ii)] A semimartingale $X$ is such that \\
			$ \mathbb{E}\big[|X_{0}|^{2}+<M>_{T}+[\text{Var}_{[0,T]}(A)]^{2}\big]<+\infty$ if and only if $||X||<+\infty$.
		\end{itemize}
	\end{theorem}
	
	Let us now introduce some kind of reflected SDE to which we are going to prove existence and uniqueness results. We begin by the following. 
	
	\begin{definition}
		Let $(L_{t})_{0\leq t\leq T}$, $(U_{t})_{0\leq t\leq T}$ be rcll adapted processes with $L\leq U$, $(H_{t})_{0\leq t\leq T}$  be a semimartingale such that $L_{0}\leq H_{0}\leq U_{0}  $ and let $\sigma$ and $a$ be two measurable random functions defined on $\Omega\times[0,T]\times\mathbb{R}$. A triplet $(X, K^{+}, K^{-})$ is said to be a solution to the SDE with two reflecting barriers $L$ and $U$, labeled later as $E(\sigma, a, L, U)$, if :
		\begin{itemize}
			\item[(i)] $(X_{t})_{0\leq t\leq T}$ is a rcll semimartingale, 
			
			\item[(ii) ] $K^{+}$ and $K^{-}$ are rcll, adapted, non negative and increasing processes, with $K_{0}^{+}=K_{0}^{-}=0$,
			\item[(iii)] \begin{equation}L_{t}\leq X_{t} \leq U_{t} \;\;\text{for}\;\;t\in [0,T]; \label{10} \end{equation}
			\item[(iv)] \begin{equation} \displaystyle\int_{0}^{T} \big(X_{t}-L_{t}\big) \; \rm{d}K^{+}_{t}=\displaystyle\int_{0}^{T} \big(U_{t}-X_{t}\big)\; \rm{d}K^{-}_{t}=0;\label{11} \end{equation}
			\item[(v)] 
			\begin{equation}
			X_{t}=H_{t}+\displaystyle\int_{0}^{t} \sigma(s,X_{s})\rm{d}B_{s}+\displaystyle\int_{0}^{t}a(s,X_{s}) \rm{d}s+K^{+}_{t}-K^{-}_{t} \; \text{for}\; t\in [0,T]. \label{12}
			\end{equation}
		\end{itemize}
	\end{definition}

	From now on we make the following assumptions.\\
	
	\begin{itemize}
		
		\item[$(\mathbf{A_{1}})$-] For fixed $x$, $\sigma(.\; , .\; ,x)$ and $a(.\;,.\;,x)$ are rcll and adapted, and there exist $\lambda>0 $ such that, $\mathbb{P}$-a.s, for all $t\in [0,T]$,
		$$\begin{array}{ccc}
		\mathbf{(A_{1-1})}&:& \forall x\in\mathbb{R}, \;\;\;|\sigma(t,x)-\sigma(t,y)|+ |a(t,x)-a(t,y)|\leq \lambda |x-y| ,\\	
		\mathbf{(A_{1-2})}&:&\forall x\in\mathbb{R} \;\;\; |\sigma(t,x)|+|a(t,x)|\leq \lambda(1+|x|). 
		\end{array}$$
		\item[$(\mathbf{A_{2}})$-] The real-valued processes $L$ and $U$ are rcll adapted processes such that:
		\begin{multline*}
		C_{L,U}:=||L||_{\mathcal{S}^{2}}+||U||_{\mathcal{S}^{2}}\\
		+ \sup_{\pi}\mathbb{E}\Big[\Big(\displaystyle\sum_{i=0}^{m-1}\big(\mathbb{E}_{\tau_{i}} (U_{\tau_{i+1}})-L_{\tau_{i}} \big)^{+}+\big(U_{\tau_{i}}-\mathbb{E}_{\tau_{i}}(L_{\tau_{i+1}}) \big)^{+}\Big)^{2}\Big]<+\infty,
		\end{multline*}
		the supremum is taken over all partitions $\pi: 0=\tau_{0}\leq ...\leq \tau_{m}=T$ of $[0,T]$.
		
	\end{itemize}

	\begin{remark}
		The integrals in (\ref{12}) are well defined. Indeed, since the process $(X_{t})_{0\leq t\leq T}$ is rcll and adapted and due to Lipschitz condition $(\mathbf{A_{1-1}})$ the processes $\big(\sigma(t,X_{t})\big)_{0\leq t\leq T}$ and $\big(a(t,X_{t})\big)_{0\leq t\leq T}$ are rcll and adapted,  for instance see \cite{4}. Furthermore, combining (\ref{10}) and the condition $(\mathbf{A_{1-2}})$, we get: 
		$$\mathbb{E}\bigg[\displaystyle\int_{0}^{T}\big(\sigma(s,X_{s})\big)^{2}ds\bigg]< +\infty \;\;\;\text{and}\;\;\; \displaystyle\int_{0}^{T} |a(s,X_{s})|\rm{d}s <+\infty \;\;\mathbb{P}\text{-a.s} .$$
		
	\end{remark}
	
	\begin{theorem}
		Under the assumptions $(\mathbf{H_{1}})$, $(\mathbf{A_{1}})$ and $(\mathbf{A_{2}})$, the reflected SDE  $E(\sigma, a, L, U)$ has at most one solution.
	\end{theorem}
	\begin{pve}
		Suppose that there exists two triplets $(X, K^{+}, K^{-})$ and  $(X^{*}, K^{+,*}, K^{-,*})$ solutions to $E(\sigma, a, L, U)$. Let $t\in [0,T]$, and put\\ $D_{t}=(K_{t}^{+}-K_{t}^{-})-(K_{t}^{+,*}-K_{t}^{-,*})$, so from (\ref{12}) $\Delta(X-X^{*})_{s}= \Delta D_{s}$. By applying It\^o's formula to the semimartingale $(X_{t}-X^{*}_{t})_{0\leq t\leq T}$ we obtain
		\begin{equation*}
		\begin{split}
		(X_{t}-X_{t}^{*})^{2}  &= 2\displaystyle\int_{0}^{t}(X_{s_{-}}-X^{*}_{s_{-}})\rm{d}(X-X^{*})_{s}\,+\displaystyle\int_{0}^{t}\big(\sigma(s,X_{s})-\sigma(s,X^{*}_{s})\big)^{2}ds \\
		& +\displaystyle\sum_{0< s\leq t}\,\big(\Delta(X-X^{*})_{s}\big)^{2}\\
		&= 2\displaystyle\int_{0}^{t} (X_{s_{-}}-X^{*}_{s_{-}})(\sigma(s,X_{s})-\sigma(s,X^{*}_{s})) \rm{d}B_{s}\\
		&+2\displaystyle\int_{0}^{t}(X_{s}-X_{s}^{*})(a(s,X_{s})-a(s,X_{s}^{*}))\rm{d}s
		+\displaystyle\int_{0}^{t}\big(\sigma(s,X_{s})-\sigma(s,X^{*}_{s})\big)^{2}\rm{d}s\\ &+2\displaystyle\int_{0}^{t}(X_{s^{-}}-X_{s^{-}}^{*}) \rm{d}D_{s}
		+\displaystyle\sum_{0< s\leq t}\big(\Delta(X-X^{*})_{s}\big)^{2}
		\end{split}
		\end{equation*}
		
		On the one hand, by adding and subtracting $2\displaystyle\int_{0}^{t}(X_{s}-X_{s}^{*})\rm{d}D_{s}$ we obtain
		\begin{equation*}
		\begin{split}
		(X_{t}-X_{t}^{*})^{2}&= 2\displaystyle\int_{0}^{t} (X_{s_{-}}-X^{*}_{s_{-}})(\sigma(s,X_{s})-\sigma(s,X^{*}_{s})) \rm{d}B_{s}\\
		&  +2\displaystyle\int_{0}^{t}(X_{s}-X_{s}^{*})(a(s,X_{s})-a(s,X_{s}^{*}))\rm{d}s
		+ \displaystyle\int_{0}^{t}\big(\sigma(s,X_{s})-\sigma(s,X^{*}_{s})\big)^{2}\rm{d}s\\
		& +2\displaystyle\int_{0}^{t}(X_{s}-X_{s}^{*}) \rm{d}D_{s}-\displaystyle\sum_{0< s\leq t}\big(\Delta(X-X^{*})_{s}\big)^{2}.
		\end{split}
		\end{equation*}
		On the other hand, using (\ref{11}) we have
		\begin{multline*}
		\displaystyle\int_{0}^{t}(X_{s}-X_{s}^{*}) \rm{d}D_{s}=-\displaystyle\int_{0}^{t}(X_{s}^{*}-L_{s}) \rm{d}K^{+}_{s}-\displaystyle\int_{0}^{t} (U_{s}-X_{s}^{*}) dK^{-}_{s}\\-\displaystyle\int_{0}^{t}(X_{s}-L_{s}) \rm{d}K^{+,*}_{s}
		-\displaystyle\int_{0}^{t} (U_{s}-X_{s}) \rm{d}K^{-,*}_{s},
		\end{multline*}
		
		which is smaller than zero by \eqref{10}. Thus
		
		\begin{equation*}
		\begin{split}
		(X_{t}-X_{t}^{*})^{2}& \leq 2\displaystyle\int_{0}^{t} (X_{s_{-}}-X^{*}_{s_{-}})(\sigma(s,X_{s})-\sigma(s,X^{*}_{s})) \rm{d}B_{s}\\
		&	+2\displaystyle\int_{0}^{t}(X_{s}-X_{s}^{*})(a(s,X_{s})-a(s,X_{s}^{*}))\rm{d}s+ \displaystyle\int_{0}^{t}\big(\sigma(s,X_{s})-\sigma(s,X^{*}_{s})\big)^{2}ds.
		\end{split}
		\end{equation*}
		Taking the expectation and using the fact that $a$ and $\sigma$ are $\lambda$-Lipschitz, we get
		$$
		\mathbb{E}\Big[(X_{t}-X_{t}^{*})^{2}\Big] \leq  (2\lambda+\lambda^{2}) \displaystyle\int_{0}^{t}\mathbb{E}\big[(X_{s}-X_{s}^{*})^{2}\big] \rm{d}s.
		$$
		
		And since plainly, 
		\begin{equation*}
		\mathbb{E}\Big[(X_{t}-X_{t}^{*})^{2}\Big]\leq 2\mathbb{E}\Big[ \sup_{0\leq s \leq T}|L_{s}|^{2}+\sup_{0\leq s \leq T}|U_{s}|^{2} \Big] \leq 2  C_{L,U}
		\end{equation*}
		we deduce by Gronwall's lemma that $\mathbb{E}\big[(X_{t}-X_{t}^{*})^{2}\big]=0$, and hence $X_{t}=X_{t}^{*}, \;\; \mathbb{P}\text{-a.s}$	and therefore $X$ and $X^{*}$ are indistinguishable since $X$ and $X^{*}$ are rcll. As in the proof of Theorem 1 we use the fact that $L$ and $U$ are completely separated to show that $K^{+}=K^{+,*}$ and $K^{-}=K^{-,*}$ This completes the proof of the theorem.
		
	\end{pve}
	
	Our second main result of this Section is the following.
	\begin{theorem}
		Assume the assumption $(\mathbf{H_{1}})$ holds.	Under $(\mathbf{A_{1}})$, $(\mathbf{A_{2}})$, the reflected SDE $E(\sigma, a, L, U)$ has a solution.
	\end{theorem}
	\begin{pve} The proof will be done in several steps.
		
		\textbf{Step 1.} As a first step, we define a Picard-type scheme. Let $t\in [0,T]$ we set
		$$X_{t}^{0}=H_{t}\;\; \text{and}\;\;
		Y_{t}^{1}=H_{t}+\displaystyle\int_{0}^{t} \sigma(s,X_{s}^{0})\rm{d}B_{s}+\displaystyle\int_{0}^{t} a(s,X_{s}^{0})\rm{d}s .$$
		By Theorem 1, there exists $(X^{1},K^{1,+},K^{1,-})$ with $K^{1,+},K^{1,-}$ are increasing and rcll processes such that:
		
		$$X_{t}^{1} =Y_{t}^{1}+K^{1,+}_{t}-K^{1,-}_{t},$$
		$$L_{t} \leq X_{t}^{1} \leq U_{t},$$ and 
		$$\displaystyle\int_{0}^{T}(X_{t}^{1}-L_{t})\; \rm{d}K^{+,1}_{t}=\displaystyle\int_{0}^{T} (U_{t}-X_{t}^{1})\; \rm{d}K^{-,1}_{t}=0 .$$
		By induction, for an integer  $n\geq 1$, we construct a triplet of rcll adapted processes $(X^{n},K^{n,+},K^{n,-})$, such that $(K^{n,+})_{t\in [0,T]}$ and $(K^{n,-})_{t\in [0,T]}$ are increasing processes, with $K_{0}^{n,+}=K_{0}^{n,-}=0 $ and
		\begin{numcases}{}
		X_{t}^{n}=H_{t}+\displaystyle\int_{0}^{t} \sigma(s,X_{s}^{n-1})\rm{d}B_{s}+\displaystyle\int_{0}^{t}a(s,X_{s}^{n-1})ds+K^{n,+}_{t}-K^{n,-}_{t} \label{13} \\
		L_{t} \leq X_{t}^{n}\leq U_{t} ; \label{14}\\
		\displaystyle\int_{0}^{T} (X_{s}^{n}-L_{s}) \rm{d}K_{s}^{n,+} =\displaystyle\int_{0}^{T} (U_{s}-X_{s}^{n}) \rm{d}K_{s}^{n,-}=0. \label{15}
		\end{numcases}
		
		Whence the Picard-type scheme we consider is well defined.\\
		
		All we have to do is to show the convergence of our scheme to the solution of $E(\sigma, a, L, U)$. It  should be noted that the difficulty arising here come from the convergence of the  sequence $K^{n}:=K^{n,+}-K^{n,-}$ to a process of bounded variation. This is why we introduce the norm $||.||$ which allows us to control the variation of $K^{n}$ by controlling $X^{n}$ 
		
		\textbf{Step 2.} The sequence of processes $(X^{n})_{n\geq 1}$ and $(K^{n})_{n\geq 1}$ converge in $\mathcal{S}^{2}$.\\
		Let $n\geq 1$. We notice first that $X_{n}\in\mathcal{S}^{2}$ due to the inequality (\ref{14}) and assumption $\mathbf{(A_{2})}$. 
		Applying It\^o's formula,  we have:
		\begin{equation*}
		\begin{split}
		(X_{t}^{n+1}-X_{t}^{n})^{2}& \leq  2\displaystyle\int_{0}^{t} (X_{s^{-}}^{n+1}-X^{n}_{s^{-}})(\sigma(s,X_{s}^{n})-\sigma(s,X^{n-1}_{s})) \rm{d}B_{s}\\ &+2\displaystyle\int_{0}^{t}(X_{s}^{n+1}-X_{s}^{n})(a(s,X_{s}^{n})-a(s,X_{s}^{n-1}))\rm{d}s\\
		& 	 +\displaystyle\int_{0}^{t}\big(\sigma(s,X_{s}^{n})-\sigma(s,X^{n-1}_{s})\big)^{2}\rm{d}s,
		\end{split}
		\end{equation*}
		which is smaller than 
		\begin{equation*}
		\begin{split}
		& 2\displaystyle\sup_{0\leq u \leq t}\Big|\displaystyle\int_{0}^{u} (X_{s^{-}}^{n+1}-X^{n}_{s^{-}})(\sigma(s,X_{s}^{n})-\sigma(s,X^{n-1}_{s})) \rm{d}B_{s}\Big|\\&+ 2\displaystyle\int_{0}^{t}\big|(X_{s}^{n+1}-X_{s}^{n})(a(s,X_{s}^{n})-a(s,X_{s}^{n-1}))\big|\rm{d}s
		+\displaystyle\int_{0}^{t}\big(\sigma(s,X_{s}^{n})-\sigma(s,X^{n-1}_{s})\big)^{2}\rm{d}s.\\
		\end{split}
		\end{equation*}
		Taking the expectation and applying Burkholder-Davis-Gundy's inequality, and the fact that $\sigma$ and $a$ are Lipschitz, we obtain
		
		$\begin{array}{lll}
		\mathbb{E}\Big[\displaystyle\sup_{0\leq s\leq t} (X_{s}^{n+1}-X_{s}^{n})^{2}\Big]&\leq&2C_{1}\lambda\mathbb{E}\Bigg[\displaystyle\sup_{0\leq s\leq t}|(X_{s}^{n+1}-X_{s}^{n})| \Bigg( \displaystyle\int_{0}^{t} \big(X_{s}^{n}-X_{s}^{n-1}\big)^{2}\rm{d}s\Bigg)^{\frac{1}{2}} \Bigg]\\
		&+& 2\lambda\mathbb{E}\Big[\displaystyle\sup_{0\leq s\leq t}|(X_{s}^{n+1}-X_{s}^{n})| \displaystyle\int_{0}^{t}\big|X_{s}^{n}-X_{s}^{n-1}\big|\rm{d}s \Big]\\
		&+& \lambda^{2}\mathbb{E}\Big[\displaystyle\int_{0}^{t} \big|X_{s}^{n}-X_{s}^{n-1}\big|^{2} \rm{d}s\Big]\\
		&\leq& \alpha \mathbb{E}\Big[\displaystyle\sup_{0\leq s\leq t} (X_{s}^{n+1}-X_{s}^{n})^{2}\Big]+\frac{C_{1}^{2}\lambda^{2}}{\alpha}\mathbb{E}\Bigg[ \displaystyle\int_{0}^{t} \big(X_{s}^{n}-X_{s}^{n-1}\big)^{2}\rm{d}s\Bigg]\\
		&+&\beta \mathbb{E}\Big[\displaystyle\sup_{0\leq s\leq t} (X_{s}^{n+1}-X_{s}^{n})^{2}\Big]+\frac{\lambda^{2}T}{\beta}\mathbb{E}\Bigg[ \displaystyle\int_{0}^{t} \big(X_{s}^{n}-X_{s}^{n-1}\big)^{2}\rm{d}s\Bigg]\\
		&+& \lambda^{2}\mathbb{E}\Bigg[ \displaystyle\int_{0}^{t} \big(X_{s}^{n}-X_{s}^{n-1}\big)^{2}\rm{d}s\Bigg],
		\end{array}$
		
		where $\alpha\;\text{and}\;\beta$ are constants strictly positive  with $\alpha+\beta<1$. So we get:
		$$\mathbb{E}\Big[\sup_{0\leq s\leq t} (X_{s}^{n+1}-X_{s}^{n})^{2}\Big]\leq M\mathbb{E}\Bigg[ \displaystyle\int_{0}^{t} \sup_{0\leq u\leq s}\big(X_{u}^{n}-X_{u}^{n-1}\big)^{2}\rm{d}s\Bigg], $$ 
		with $M:=\big(\frac{C_{1}^{2}\lambda^{2}}{\alpha}+\frac{T\lambda^{2}}{\beta}+\lambda^{2}\big)(1-\alpha-\beta)^{-1}$.
		
		We denote $h_{n}(t):=\mathbb{E}\Big[\displaystyle\sup_{0\leq u\leq t} (X_{u}^{n}-X_{u}^{n-1})^{2}\Big]$, and we put $h_{0}(t):=\mathbb{E}\Big[\displaystyle\sup_{0\leq u\leq t} (X_{u}^{0})^{2}\Big]$
		
		Since $h_{n+1}(t)\leq M\displaystyle\int_{0}^{t}h_{n}(s)\rm{d}s$, an induction argument shows that $h_{n}(t)\leq (h_{0}(t))M^{n} \frac{t^{n}}{n!}$. and since $\displaystyle\sum_{n\geq 0} \big((h_{0}(t))M^{n} \frac{t^{n}}{n!}\big)^{\frac{1}{2}}< +\infty$ then $$ \displaystyle\sum_{n\geq 0} || (X^{n+1}-X^{n})||_{\mathcal{S}^{2}}<+\infty,$$ 
		implying that $(X^{n})_{n\geq 1}$is a Cauchy sequence in the complete space $\mathcal{S}^{2}$ and thus converges to a process $X\in\mathcal{S}^{2}$. 
		
		To show the convergence of the sequence $(K^{n})_{n\geq 1}$, we apply Burkholder-Davis-Gundy's inequality, the fact that $\sigma$ and $a$ are Lipschitz and the convergence of $(X^{n})_{n\geq 1}$ in $\mathcal{S}^{2}$ to get  
		\begin{multline*}
		\lim\limits_{n\rightarrow +\infty} \int_{0}^{t}\sigma(s,X_{s}^{n}) \rm{d}B_{s}\overset{\mathcal{S}^{2}}{=}\lim\limits_{n\rightarrow +\infty}\int_{0}^{t}\sigma(s,X_{s}) \rm{d}B_{s}\;
		\;\text{and}\\ \;\lim\limits_{n\rightarrow +\infty} \int_{0}^{t}a(s,X_{s}^{n}) \rm{d}s\overset{\mathcal{S}^{2}}{=}\lim\limits_{n\rightarrow +\infty}\int_{0}^{t}a(s,X_{s}) \rm{d}s.
		\end{multline*}
		
		The equation (\ref{13}) implies $(K^{n})_{n\geq 1}$ converges in $\mathcal{S}^{2}$ to a process $K\in\mathcal{S}^{2}$.
		\\
		
		\textbf{Step 3.}  We show that the process $K$ is a process of bounded variation.\\
		Since $K^{n}$ converges in $\mathcal{S}^{2}$ to a process $K\in\mathcal{S}^{2}$, there exists a subsequence $(n_{k})_{k\geq 1}$ such that: 
		$$\lim\limits_{k\rightarrow +\infty}\sup_{0\leq t\leq T} \mid K^{n_{k}}_{t}-K_{t}\mid = 0\;\; \mathbb{P}\text{- a.s} .$$
		
		On the other hand, using the fact that  $L\leq X^{n_{k}} \leq U$,  and the supremum is again taken over all partitions $\pi:0=\tau_{0}\leq...\leq\tau_{m}=T$  we obtain
		\begin{multline*}
		\displaystyle\sup_{\pi}\mathbb{E}\Bigg[\Big( \displaystyle\sum_{i=0}^{m-1}|\mathbb{E}_{\tau_{i}}(X_{\tau_{i+1}}^{n_{k}})-X_{\tau_{i}}^{n_{k}}|\Big)^{2}\Bigg] \leq\\ \sup_{\pi}\mathbb{E}\Bigg[\Big(\displaystyle\sum_{i=0}^{m-1}\big(\mathbb{E}_{\tau_{i}} (U_{\tau_{i+1}})-L_{\tau_{i}} \big)^{+}+\big(U_{\tau_{i}}-\mathbb{E}_{\tau_{i}}(L_{\tau_{i+1}})\big)^{+}\Big)^{2}\Bigg],
		\end{multline*}

		which is finite from the assumption $\mathbf{(A_{2})}$.\\
		
		Applying the inequality in the point  $(i)$ of Theorem 2, we get:
		$$\mathbb{E}((K_{T}^{n_{k},+}+K_{T}^{n_{k},+})^{2})=\mathbb{E}([\text{Var}_{[0,T]}(K^{n_{k}})]^{2}) \leq C ||X^{n_{k}}||\leq C C_{L,U}<+\infty.$$
		
		Hence  by the section theorem (see, Theorem 4.12 of [\cite{7}, p 116 ]), we get:
		$$(K_{T}^{n_{k},+}+K_{T}^{n_{k},-})\leq (C C_{L,U})^{\frac{1}{2}} \;\;\; \mathbb{P}\text{-a.e}.$$
		Due to Helly's selection theorem, we conclude that $K$ is a process of bounded variation. So there exist an increasing rcll processes $K^{+}$ and $K^{-}$ such that $K=K^{+}-K^{-}$ with $dK^{+}$ and $dK^{-}$ has disjoint support.
		
		\textbf{Step 4.} We show that: $$\displaystyle\int_{0}^{T} (X_{s}-L_{s}) \rm{d}K_{s}^{+}=\displaystyle\int_{0}^{T} (U_{s}-X_{s}) \rm{d}K_{s}^{-}=0.$$
		As in \cite{11} we will show that  $K$ increase if and only if $K=L$ and decrease if and only if $K=U$. \\If $t$ is a point of increase of $K$, then there exist $\epsilon>0$ such that  $$K_{t^{-}}< K_{s}\;\; \text{for}\;\; s\in ]t,t+\epsilon].$$
		Since $K^{n_{k}}$ converges uniformly, we obtain $$\lim\limits_{k\rightarrow +\infty} K^{n_{k}}_{t^{-}}=K_{t^{-}} \;\; \mathbb{P}\text{-a.s}.$$
		So there exists $N$ such that for $k\geq N $,  $$K^{n_{k}}_{t^{-}}< K^{n_{k}}_{s}\;\; \text{for}\;\; s\in]t,t+\epsilon],$$ 
		Therefore the point $t$ is a point of increase of $K^{n_{k}}$, hence $X^{n_{k}}_{t}=L_{t}$, and then 
		$$X_{t} = L_{t} \;\; \mathbb{P}\text{- a.s} .$$
		Similarly we show that  the support of $dK^{-}$ is included in the set $\{X=U\}$.\\
		In conclusion combining all the steps above we have shown that the Picard-type scheme converges to the solution of $E(\sigma, a, L, U)$. That completes the proof.
		
	\end{pve}

	\section{Appendix}
	Here we give a slight generalization of the result in [\cite{8},p 270-271 ], for the case of rcll adapted processes. We will show  the existence and uniqueness of Skorokhod problem under the assumption $\mathbf{(H_{0})}$.
	\begin{proposition}\label{proApp}
		Let $(Y_{t})_{t\geq 0}$, $(L_{t})_{t\geq 0}$ and $(U_{t})_{t\geq 0}$ be rcll adapted processes such that $L$ and $U$ satisfies $\mathbf{(H_{0})}$. Then, the Skorokhod problem $SPR(Y, L, U)$ has a unique solution.
	\end{proposition}
	
	\begin{pve}
		We construct $K^{+}$ and $K^{-}$ by alternating between the two one-sided reflection operators corresponding to downward reflection at $L$ and upward reflection at $U$.	We assume that $Y$ hits the lower boundary first, we set:
		$$X_{t}^{0,0}=Y_{t},\;\; \phi_{t}^{1}=\sup_{s\leq t}(L_{s}-X_{s}^{0,0}) \lor 0 \;\; \text{and}\; \;\tilde{X_{t}}^{1,0}=Y_{t}+\phi_{t}^{1}.$$  
		Let $T_{1}=\inf\{ t\geq 0 \;:\; Y_{t}\leq L_{t}\}\;\text{and}\; S_{1}=\inf\{ t\geq T_{1}\;:\; \tilde{X}_{t}^{1,0}\geq U_{t}\}$
		$ \text{with the convention}\;\inf(\emptyset)=\infty $.
		
		We modify $\tilde{X}^{1,0}$ as :
		$$X _{t}^{1,0}:=Y_{t}+\phi_{t\land S_{1}}^{1}. $$
		We set:
		$$\psi_{t}^{1}=\sup_{s\leq t}(X_{s}^{1,0}-U_{s}) \lor 0 \;\; \text{and}\;\; \tilde{X}_{t}^{1,1}=Y_{t}+\phi_{t\land S_{1}}^{1} - \psi_{t}^{1}.$$

		We set $T_{2}= \inf\{ t\geq S_{1} \;:\; \tilde{X}_{t}^{1,1}\leq L_{t}\} \;\; \text{and} \;\inf(\emptyset)=\infty$, and we modify $\tilde{X}_{t}^{1,1}$ as :
		$$X_{t}^{1,1}=Y_{t}+\phi_{t\land S_{1}}^{1} -\psi_{t\land T_{2}}^{1}$$
		
		The triplet $(X^{1,1}, \phi_{.\land S_{1}}^{1}, \psi_{.\land L_{2}}^{1})$ satisfy \eqref{1}, \eqref{2}, and \eqref{3} on $[0,T_{2}[.$ 
		
		For $n\geq 1$ we define by induction:
		$$\phi_{t}^{n}=\sup_{s\leq t}(L_{s}-X_{s}^{n-1,n-1}) \lor 0\,\; \text{and}\;\;\psi_{t}^{n}=\sup_{s\leq t}(X_{s}^{n,n-1}-U_{s}) \lor 0 ,$$
		$$S_{n}=\inf\{ t\geq T_{n}  \;:\; \tilde{X}_{t}^{n,n-1}\geq U_{t}\},\; T_{n+1}=\inf\{ t\geq S_{n} \;:\; \tilde{X}^{n,n}_{t}\leq L_{t}\},\;  \inf(\emptyset)=\infty $$
		
		and
		$$X_{t}^{n+1,n+1}=Y_{t}+\sum_{k=1}^{n}\phi^{k}_{t\land S_{k}}- \sum_{k=1}^{n}\psi^{k}_{t\land T_{k+1}}. $$
		The triplet $\big(X^{n,n}, \sum_{k=1}^{n-1}\phi^{k}_{.\land S_{k}}, \sum_{k=1}^{n-1}\psi^{k}_{.\land T_{k+1}}\big) $ satisfy \eqref{1},\eqref{2}, and \eqref{3} on $[0,T_{n+1}[$.
		
		We denote by $T_{\infty}:=\lim\limits_{n\rightarrow +\infty} T_{n}$ and we set :
		$$X_{t}=Y_{t}+K_{t}^{+}-K_{t}^{-}, \;\; \text{for} \;\; t< T_{\infty}$$
		with 
		$$ K_{t}^{+}= \sum_{k=1}^{+\infty}\phi^{k}_{t\land S_{k}} \;\; \text{and}\;\;K_{t}^{-}=\sum_{k=1}^{+\infty}\psi^{k}_{t\land T_{k+1}}.$$
		In order to show that $(X,K^{+},K^{-})$ is the solution of $SPR(Y, L, U)$, we have to prove that: $$ T_{\infty}=+\infty.$$
		For an arbitrary $t$. The total number of such crossing of $X^{n,n}$ from the lower boundary to the upper boundary in $[0,t]$ is finite that we denote by $N$. We define:
		$$\tau_{i}: \text{ the}\; i^{th}\; \text{time}\; X^{n,n}\; \text{hits the lower barriers } L \;\text{on}\; [0,t]\;\; i\leq N$$
		$$r_{i}: \text{ the}\; i^{th}\; \text{time}\; X^{n,n}\; \text{hits the upper barrier}\; U \text{on}\; [0,t]\;\; i\leq N$$
		
		We have $$|X_{\tau_{i}}^{n,n}-X_{r_{i}}^{n,n}|=|L_{\tau_{i}}-U_{r_{i}}|\geq\inf_{s\leq t}(U_{s}-L_{s})>\epsilon.$$
		
		Since $\inf_{s\leq t}(U_{s}-L_{s})>\epsilon$, by Lemma 1 [(\cite{1}), p, 122], there exists $\delta>0$ such that $r_{i}-\tau_{i}>\delta,$ for all $i\leq N$ and then we get $N\delta\leq t$. Consequently $N$ is bounded by $\frac{t}{\delta}$
		and hence $t\leq T_{\infty}$. Since $t$ is arbitrary we obtain $T_{\infty}=+\infty$.
		And finally 
		$$ K_{t}^{+}= \sum_{k=1}^{+\infty}\phi^{k}_{t\land S_{k}}<\infty \;\; \text{and}\;\; K_{t}^{-}= \sum_{k=1}^{+\infty}\psi^{k}_{t\land T_{k+1}}<\infty.$$
		and the triplet $(X, K^{+}, K^{-})$ satisfy (\ref{1}), (\ref{2}), and (\ref{3}) for $t\in [0,+\infty[$.
		
		The uniqueness may be obtained in a similar way as in the proof of Theorem 1. 
	\end{pve}

\end{document}